\documentclass[11pt]{article}
\usepackage{amsmath,amssymb,theorem}
\usepackage{graphicx, enumerate}
\usepackage{subfigure}
\usepackage{psfrag}
\usepackage{placeins}
\usepackage{wrapfig}
\usepackage{booktabs}
\usepackage{chngcntr}
\counterwithin{figure}{section}
\numberwithin{figure}{section}
\counterwithin{table}{section}
\usepackage{hyperref}
\hypersetup{
  colorlinks,
  citecolor=black,
  filecolor=black,
  linkcolor=black,
  urlcolor=black 
}

\newcommand{\mc}{\mathcal}

\newtheorem{thm}{Theorem}[section]
\newtheorem{conj}[thm]{Conjecture}
\newtheorem{prop}[thm]{Proposition}
\theorembodyfont{\rmfamily}
\def\pf{\bigskip\noindent {\bf Proof.}~~}

\def\dfn#1{{\sl #1}}
\def\es{\emptyset}
\def\less{\backslash}

\newcounter{counter}

\def\proofsquare{  \bigskip\hfill\vrule height3pt width6pt depth2pt}

\begin{document}
\title{A conjecture on Gallai-Ramsey numbers of even cycles and paths}
\author{Zi-Xia Song and Jingmei Zhang\thanks{Email addresses: Zixia.Song@ucf.edu (Z-X. Song) and jmzhang@knights.ucf.edu (J. Zhang).} \\
 Department  of Mathematics\\
 University of Central Florida \\  
 Orlando, FL 32816, USA \\
}
\date{}
\maketitle

\begin{abstract}
A {\it Gallai coloring} is a coloring of the edges of a complete graph without rainbow triangles, and a {\it Gallai $k$-coloring} is a Gallai coloring that uses at most $k$ colors.  Given   an integer $k\ge1$ and graphs $H_1, \ldots, H_k$, the      Gallai-Ramsey number $GR(H_1,  \ldots, H_k)$ is the least integer $n$ such that every Gallai $k$-coloring of the complete graph $K_n$   contains a monochromatic copy of $H_i$ in color $i$ for some $i \in \{1,2, \ldots, k\}$.  When $H = H_1 = \cdots = H_k$, we simply write $GR_k(H)$.   
We study   Gallai-Ramsey numbers of  even cycles    and paths.   For all $n\ge3$ and $k\ge2$, let  $G_i=P_{2i+3}$ be  a path on $2i+3$ vertices for all $i\in\{0,1, \ldots, n-2\}$ and $G_{n-1}\in\{C_{2n}, P_{2n+1}\}$. Let   $ i_j\in\{0,1,\ldots, n-1 \}$  for all $j\in\{1,2, \ldots, k\}$ with 
 $  i_1\ge i_2\ge\cdots\ge i_k $. 
The first author recently conjectured that $ GR(G_{i_1}, G_{i_2}, \ldots, G_{i_k}) = |G_{i_1}|+\sum_{j=2}^k i_j$. The truth of this conjecture implies that  $GR_k(C_{2n})=GR_k(P_{2n})=(n-1)k+n+1$ for all $n\ge3$ and $k\ge1$,  and $GR_k(P_{2n+1})=(n-1)k+n+2$ for all $n\ge1$ and $k\ge1$.  In this paper, we prove   that  the aforementioned  conjecture holds for  $n\in\{3,4\}$ and    all $k\ge2$. Our proof relies only on Gallai's result and the classical Ramsey numbers $R(H_1, H_2)$, where $H_1, H_2\in\{C_8, C_6, P_7, P_5, P_3\}$. 
We believe the recoloring method we developed here will  be very useful for solving subsequent cases, and perhaps the conjecture. 
\end{abstract}
\section{Introduction}

\baselineskip 16pt

 In this paper we consider graphs that are finite, simple and undirected.    Given a graph $G$ and a set $A\subseteq V(G)$,  we use   $|G|$    to denote  the  number
of vertices    of $G$, and  $G[A]$ to denote the  subgraph of $G$ obtained from $G$ by deleting all vertices in $V(G)\less A$.  A graph $H$ is an \dfn{induced subgraph} of $G$ if $H=G[A]$ for some $A\subseteq V(G)$.  We use $P_n$,  $C_n$ and $K_n$ to denote the path,    cycle and  complete graph  on $n$ vertices, respectively.  
For any positive integer $k$, we write  $[k]$ for the set $\{1,2, \ldots, k\}$. We use the convention   ``$A:=$'' to mean that $A$ is defined to be the right-hand side of the relation.

 Given an integer $k \ge 1$ and graphs $H_1,  \ldots, H_k$, the classical Ramsey number $R(H_1,   \ldots, H_k)$   is  the least    integer $n$ such that every $k$-coloring of  the edges of  $K_n$  contains  a monochromatic copy of  $H_i$ in color $i$ for some $i \in [k]$.  Ramsey numbers are notoriously difficult to compute in general. In this paper, we  study Ramsey numbers of graphs in Gallai colorings, where a \dfn{Gallai coloring} is a coloring of the edges of a complete graph without rainbow triangles (that is, a triangle with all its edges colored differently). Gallai colorings naturally arise in several areas including: information theory~\cite{KG}; the study of partially ordered sets, as in Gallai's original paper~\cite{Gallai} (his result   was restated in \cite{Gy} in the terminology of graphs); and the study of perfect graphs~\cite{CEL}. There are now a variety of papers  which consider Ramsey-type problems in Gallai colorings (see, e.g., \cite{C9C11, C13C15, DylanSong,chgr, c5c6,GS, exponential, Hall}).   These works mainly focus on finding various monochromatic subgraphs in such colorings. More information on this topic  can be found in~\cite{FGP, FMO}.  
 
A \dfn{Gallai $k$-coloring} is a Gallai coloring that uses at most $k$ colors. 
 Given an integer $k \ge 1$ and graphs $H_1,  \ldots, H_k$, the   \dfn{Gallai-Ramsey number} $GR(H_1,  \ldots, H_k)$ is the least integer $n$ such that every Gallai $k$-coloring of $K_n$   contains a monochromatic copy of $H_i$ in color $i$ for some $i \in [k]$. When $H = H_1 = \dots = H_k$, we simply write $GR_k(H)$ and  $R_k(H)$.    Clearly, $GR_k(H) \leq R_k(H)$ for all $k\ge1$ and $GR(H_1, H_2) = R(H_1, H_2)$.    In 2010, 
Gy\'{a}rf\'{a}s,   S\'{a}rk\"{o}zy,  Seb\H{o} and   Selkow~\cite{exponential} proved   the general behavior of $GR_k(H)$.

\begin{thm} [\cite{exponential}]
Let $H$ be a fixed graph  with no isolated vertices 
 and let $k\ge1$ be an integer. Then
$GR_k (H)$ is exponential in $k$ if  $H$ is not bipartite,    linear in $k$ if $H$ is bipartite but  not a star, and constant (does not depend on $k$) when $H$ is a star.			
\end{thm}

It turns out that for some graphs $H$ (e.g., when $H=C_3$),  $GR_k(H)$ behaves nicely, while the order of magnitude  of $R_k(H)$ seems hopelessly difficult to determine.  It is worth noting that  finding exact values of $GR_k (H)$ is  far from trivial, even when $|H|$ is small.
We will utilize the following important structural result of Gallai~\cite{Gallai} on Gallai colorings of complete graphs. 
\begin{thm}[\cite{Gallai}]\label{Gallai}
	For any Gallai coloring $c$ of a complete graph $G$ with $|G|\ge2$, $V(G)$ can be partitioned into nonempty sets  $V_1,  \dots, V_p$ with $p>1$ so that    at most two colors are used on the edges in $E(G)\less (E(G[V_1])\cup \cdots\cup  E(G[V_p]))$ and only one color is used on the edges between any fixed pair $(V_i, V_j)$ under $c$. \end{thm}

The partition given in Theorem~\ref{Gallai} is  a \dfn{Gallai partition} of  the complete graph $G$ under  $c$.  Given a Gallai partition $V_1, \dots, V_p$ of the complete graph $G$ under $c$, let $v_i\in V_i$ for all $i\in[p]$ and let $\mathcal{R}:=G[\{v_1,  \dots, v_p\}]$. Then $\mathcal{R}$ is  the \dfn{reduced graph} of $G$ corresponding to the given Gallai partition under $c$. Clearly,  $\mathcal{R}$ is isomorphic to $K_p$.  It is worth noting that  $\mathcal{R}$ does not depend on the choice of $v_1, \ldots, v_p$ because $\mathcal{R}$ can be   obtained by first contracting each part $V_i$ into a single vertex, say $v_i$,   and then coloring every edge  $v_iv_j$ by the color used on the edges between $V_i$ and $V_j$ under $c$.  
By Theorem~\ref{Gallai},  all the edges in $\mathcal{R}$ are colored by at most two colors under $c$.  One can see that any monochromatic copy of $H$ in $\mathcal{R}$ under $c$ will result in a monochromatic copy of $H$ in $G$ under $c$. It is not  surprising that  Gallai-Ramsey numbers $GR_k(H)$ are closely related to  the classical Ramsey numbers $R_2(H)$.  Recently,  Fox,  Grinshpun and  Pach~\cite{FGP} posed the following  conjecture on $GR_k(H)$ when $H$ is a complete graph. 

\begin{conj}[\cite{FGP}]\label{Fox} For all    $t\ge3$ and $k\ge1$,
\[
GR_k(K_t) = \begin{cases}
			(R_2(K_t)-1)^{k/2} + 1 & \text{if } k \text{ is even} \\
			(t-1)  (R_2(K_t)-1)^{(k-1)/2} + 1 & \text{if } k \text{ is odd.}
			\end{cases}
\]
\end{conj}

Recall that if $n< R_k(K_3)$, then there is a $k$-coloring $c$ of the edges of $K_n$ such that edges of every triangle  in $K_n$ are colored by at least two colors under $c$.  A question of T. A. Brown (see \cite{chgr}) asked: What is the largest number $f(k)$ of vertices of a complete graph can have such that it is possible to $k$-color its edges so that edges of every  triangle are colored by   exactly two colors? Chung and Graham~\cite{chgr}  answered this question in  1983.

\begin{thm}[\cite{chgr}]\label{C3} For all  $k \ge 1$, 
	$f(k) = \begin{cases}
			5^{k/2}   & \text{if } k \text{ is even} \\
			2 \cdot  5^{(k-1)/2}   & \text{if } k \text{ is odd.}
			\end{cases}$
\end{thm}

Clearly, $GR_k(K_3)=f(k)+1$. By Theorem~\ref{C3},    Conjecture~\ref{Fox} holds for $t=3$.  The proof of Theorem~\ref{C3} does not rely on Theorem~\ref{Gallai}.  A simpler proof of this case using Theorem~\ref{Gallai}  can be found in~\cite{exponential}.    The next open case, when $t=4$,  was recently settled in~\cite{K4}.     Gallai-Ramsey number  of $H$, where $H\in\{C_4,  P_5, C_6, P_6\}$,  has also been studied, as well as  general  upper bounds for $GR_k (P_n)$ and $GR_k (C_{n})$ that  were first studied in \cite{FGJM, c5c6} and later   improved  in  \cite{Hall}.    Gregory~\cite{c8}   proved in his  thesis   that $GR_k (C_8)=3k+5$,  but the proof was incomplete.   We list some results in \cite{FGJM, c5c6, Hall} below.

\begin{thm}[\cite{FGJM}]\label{gallai-path}
For all   $k \ge 1$,  
\begin{enumerate}[\rm(a)]
\item $GR_k(P_n)=\lfloor \frac{n-2}{2}\rfloor k+ \lceil \frac n 2 \rceil +1$ for     $ n\in\{3, 4,5,6\}$.
 
\item  $GR_k(C_4) = k + 4$.
\end{enumerate}
\end{thm}

\begin{thm}[\cite{c5c6}]\label{C5C6} For all $k \ge 1$, 
$GR_k(C_5) = 2^{k+1} + 1$ and $GR_k(C_6) = 2k + 4$.
\end{thm}
 
\begin{thm}[\cite{Hall}]\label{Hall-path}
For all  $n \ge 3$ and $k \ge 1$, $$GR_k (C_{2n}) \le (n-1)k+3n  \text{ and } 
  GR_k (P_n) \le  \left\lfloor \frac{n-2}{2} \right\rfloor k +3 \left\lfloor \frac{n}{2} \right\rfloor.$$ 
\end{thm}

More recently,  Gallai-Ramsey numbers of  odd cycles on at most $15$ vertices    have been completely settled by Bruce and Song~\cite{DylanSong} for $C_7$, Bosse and Song~\cite{C9C11} for $C_9$ and $C_{11}$,  and Bosse, Song and Zhang~\cite{C13C15} for $C_{13}$ and $C_{15}$.    Very recently, the exact values of $GR_k(C_{2n+1})$ for $n \ge 8$ has been solved by Zhang,   Song and Chen~\cite{Fangfang}.  We summarize these results  below.

\begin{thm}[\cite{ C9C11, C13C15,DylanSong}] For $n\in\{3,4,5,6,7\}$ and all $k \ge 1$, 
$GR_k(C_{2n+1}) = n\cdot 2^{k} + 1$.
\end{thm}

In this paper, we study   Gallai-Ramsey numbers of   even cycles and paths. Note that  $GR_k(H)=|H|$ for any graph $H$ when $k=1$.   
 For   all  $n \ge 3$ and $k\ge2$, let $G_{n-1}\in \{C_{2n}, P_{2n+1}\}$,   $G_i :=P_{2i+3}$ for all  $i \in \{0, 1, \ldots, n-2\}$. We want to determine  the exact values of $GR (G_{i_1}, \ldots, G_{i_k})$, where    $ i_j\in\{0,1,\ldots, n-1 \}$  for all $j\in[k]$. By reordering colors if necessary, we   assume that $i_1\ge i_2\ge \cdots \ge i_k$.    The construction for establishing a lower bound for  $GR (G_{i_1},  \ldots, G_{i_k})$ for all $n\ge3$ and $k\ge2$ is similar to the construction given by   Erd\H{o}s, Faudree, Rousseau and Schelp  in 1976 (see Section 2  in \cite{EFRS}) for classical Ramsey numbers of even cycles and paths. We recall their construction  in the proof of Proposition~\ref{lower}.  
We list below   the results on $2$-colored Ramsey numbers of even cycles and paths that will be used  in the proofs of Proposition~\ref{lower} and  Theorem~\ref{main}. 

\begin{thm}[\cite{Rosta}]\label{cycles}
For all $n \ge 3$, $R_2(C_{2n})=3n-1$. 
\end{thm}

\begin{thm}[\cite{FLPS}]\label{path-cycle}
For all integers $n, m$ satisfying $2n \ge m \ge 3$, $R(P_m, C_{2n})=2n+ \lfloor \frac m2 \rfloor -1$. 
\end{thm} 

\begin{thm}[\cite{GeGy}]\label{paths}
For all integers $n, m$ satisfying $n \ge m \ge 2$, $R(P_m, P_n)=n+ \lfloor \frac m2 \rfloor -1$. 
\end{thm}

\begin{prop} \label{lower}
For   all  $n \ge 3$ and  $k\ge2$, 
$$GR(G_{i_1},  \ldots, G_{i_k}) \ge |G_{i_1}|+\sum_{j=2}^k i_j,$$
where $n-1\ge i_1\ge \cdots \ge i_k\ge0$.
\end{prop}

\pf    By Theorem~\ref{cycles},  Theorem~\ref{path-cycle} and Theorem~\ref{paths}, the statement  is   true when $k=2$. So we may assume that $k\ge3$.  To show that 
$GR (G_{i_1},  \ldots, G_{i_k}) \ge |G_{i_1}|+\sum_{j=2}^k i_j$, we recall the construction given in \cite{EFRS}.  Let $G$ be a complete graph on $ (|G_{i_1}|-1)+\sum_{j=2}^k i_j$ vertices.   Let $V_1, \ldots, V_k$ be a partition of $V(G)$ such that $|V_1|=|G_{i_1}|-1$ and $|V_j|= i_j$ for all $j\in\{2,3, \ldots, k\}$. 
Let $c$ be a   $k$-edge-coloring of  $G$   by first coloring all the edges of  $G[V_j]$ by  color $j$ for all $j\in[k]$, and  then coloring all the edges between $V_{j+1}$ and $\bigcup_{\ell=1}^jV_\ell $ by   color $j+1$ for all $j\in [k-1]$. Then   $G$ contains neither a rainbow triangle nor a monochromatic copy of $G_{i_j}$ in  color $j$ for all $j\in[k]$ under $c$. Hence, $ GR (G_{i_1},   \ldots, G_{i_k}) \ge |G|+1= |G_{i_1}|+\sum_{j=2}^k i_j $, as desired. \proofsquare

Motivated by the work developed  in \cite{c8}, the first author recently conjectured that the lower bound established in Proposition~\ref{lower} is also the desired  upper bound for $GR(G_{i_1},  \ldots, G_{i_k})$ for all $n\ge3$ and $k\ge2$. 
We state it below (note that Conjecture~\ref{Song} was first mentioned at the   49th Southeastern International Conference on Combinatorics, Graph Theory \& Computing,  Florida Atlantic University, Boca Raton, FL, March 5-9, 2018). 
\begin{conj}\label{Song}
For   all  $n \ge 3$ and $k\ge2$, 
$$GR(G_{i_1},  \ldots, G_{i_k}) =  |G_{i_1}|+\sum_{j=2}^k i_j,$$
where $n-1\ge i_1\ge  \cdots \ge i_k\ge0$.
\end{conj}

  Clearly, $GR_k(C_{2n})\ge GR_k(P_{2n})$ and $GR_k(C_{2n})\ge GR_k(M_{n})$, where   $M_n$ denotes  a matching of size $n$.  It is worth noting that by letting $i_1=\cdots=i_k=n-1$ and $G_{i_1}=C_{2n}$, the construction given in the proof of Proposition~\ref{lower} yields that $(n-1)k+n+1\le GR_k(P_{2n})$ and $(n-1)k+n+1\le GR_k(M_{n})$ for all $n\ge3$ and $k\ge1$ (the authors would like to thank Joseph Briggs, a Ph.D. student at the Carnegie-Mellon University, for pointing this out for  $M_n$,  at the   49th Southeastern International Conference on Combinatorics, Graph Theory \& Computing,  Florida Atlantic University, Boca Raton, FL, March 5-9, 2018).  
The truth of Conjecture~\ref{Song} implies that  $GR_k(C_{2n})=GR_k(P_{2n})=GR_k(M_{n})=(n-1)k+n+1$ for all $n\ge3$ and $k\ge1$ and $GR_k(P_{2n+1})=(n-1)k+n+2$ for all $n\ge1$ and $k\ge1$.  As observed in \cite{Hall},  to completely solve Conjecture~\ref{Song}, one only needs to consider the case  $G_{n-1}=C_{2n}$. We prove this    in Proposition~\ref{C2n to P2n+1}. The proof of Proposition~\ref{C2n to P2n+1} is similar to the proof of Theorem 7 given in \cite{Hall}. We include a proof here for completeness.  

\begin{prop}\label{C2n to P2n+1}
For   all  $n \ge 3$ and $k\ge2$, if Conjecture~\ref{Song} holds for     $G_{n-1}=C_{2n}$, then 
 it also holds for    $G_{n-1}=P_{2n+1}$. 
\end{prop}

\pf  By the assumed truth of Conjecture~\ref{Song}   for     $G_{n-1}=C_{2n}$,  we may assume that $G_{i_1}=P_{2n+1}$. Then  $i_1=n-1$.  We may further assume that    $n-1= i_1  = \cdots = i_t > i_{t+1} \ge \cdots \ge i_k$, where $t \in [k]$.  By Proposition~\ref{lower}, $ GR (G_{i_1},  \ldots, G_{i_k}) \ge (2n+1)+\sum_{j=2}^k i_j=2+n+t(n-1)+\sum_{j=t+1}^k i_j$.  We next show that 
$GR(G_{i_1},  \ldots, G_{i_k}) \le  2+n+t(n-1)+\sum_{j=t+1}^k i_j$. \medskip

 Let $G$ be a complete graph on $2+n+t(n-1)+\sum_{j=t+1}^k i_j$ vertices and let $c : E(G) \rightarrow [k]$ be any Gallai coloring of $G$. Suppose $G$ does not contain a monochromatic copy of $G_{i_j}$   in color $j$  for all $j \in [k]$.     By the assumed truth of Conjecture~\ref{Song}   for    $G_{n-1}=C_{2n}$, $GR (C_{2n}, \ldots, C_{2n}, G_{i_{t+1}}, \ldots, G_{i_k}) =  2n+(t-1)(n-1)+\sum_{j=t+1}^k i_j=1+n+t(n-1)+\sum_{j=t+1}^k i_j$.   Thus $G$ must contain a monochromatic  copy of    $H:=C_{2n}$ in some color $\ell \in [t]$ under $c$.  We may assume that $\ell = 1$.  Then for every vertex $u \in V(G) \less V(H)$, all the edges between $u$ and $V(H)$   must be   colored by exactly one color $j$ for some   $j\in \{2, \ldots, k\}$, because $G$   contains neither a  rainbow triangle nor a monochromatic copy of $P_{2n+1}$ in color $1$ under $c$. Thus, $V(G) \less V(H)$ can be partitioned into $V_2, V_3, \ldots, V_{k }$ such that all the edges between $V_j$ and $V(H)$ are colored   by color $j$ for all $j\in \{2, \ldots, k\}$. It follows that for all $j\in \{2, \ldots, k\}$,    $|V_j| \le i_j$, because $G$ does not contain a monochromatic copy of $G_{i_j}$ in color $j$. But then $|G|=|H|+\sum_{j=2}^{k }|V_j| \le 2n+\sum_{j=2}^{k } i_j=1+n+t(n-1)+\sum_{j=t+1}^k i_j$, contrary to $|G|=2+n+t(n-1)+\sum_{j=t+1}^k i_j$.\proofsquare

In this paper, we  prove that Conjecture~\ref{Song} is true for $n\in\{3,4\}$ and all $k\ge1$. 

\begin{thm}\label{main}
For     $n \in\{3,4\}$ and all $k\ge2$,  let  $G_i=P_{2i+3}$ for all  $i \in \{0, 1, \ldots, n-2\}$, $G_{n-1}=C_{2n}$, and     $ i_j\in\{0,1,\ldots, n-1 \}$  for all $j\in[k]$ with $i_1\ge i_2\ge \cdots \ge i_k$. Then 
$$GR(G_{i_1},  \ldots, G_{i_k}) = |G_{i_1}|+\sum_{j=2}^k i_j.$$
\end{thm}

 Theorem~\ref{main} strengthens  the results listed in Theorem~\ref{gallai-path},  Theorem~\ref{C5C6} and $GR_k(C_8)=3k+5$ given in \cite{c8}.  Our proof   relies only on Theorem~\ref{Gallai} and   Ramsey numbers $R(H_1, H_2)$, where $H_1, H_2\in\{C_8, C_6, P_7, P_5, P_3\}$.     Theorem~\ref{main},  together with Proposition~\ref{C2n to P2n+1}, implies that  $GR_k(C_{2n})=GR_k(P_{2n})=GR_k(M_{n})=(n-1)k+n+1$ for $n\in\{3,4\}$ and all $k\ge1$,  and $GR_k(P_{2n+1})=(n-1)k+n+2$ for $n\in\{1,2, 3,4\}$ and all $k\ge1$.   
Hence, Theorem~\ref{main}    yields a new and simpler proof of 
the known results on  Gallai-Ramsey numbers of  $C_8$, $C_6$ and $P_n$ with  $n\le 7$. As mentioned earlier,   the proof of $GR_k(C_8)=3k+5$  given in \cite{c8}   was incomplete.   
We prove Theorem~\ref{main} in Section~\ref{C8}.  In our completely new strategy, we developed an extremely useful recoloring method (in the proof of Claim~\ref{e:A1} in Section~\ref{C8}) which we believe will assist in solving other cases, and possibly the conjecture.
This method, together with new ideas, has been applied in \cite{SZ} to prove that Conjecture~\ref{Song}  is true for  $n\in\{5,6\}$ and all $k\ge2$. Note  that the method we developed here for even cycles and paths  is very different from the method  for odd cycles developed in \cite{C9C11, C13C15,DylanSong}.  \medskip

\section{Proof of Theorem~\ref{main}}\label{C8}

We are ready to prove Theorem~\ref{main}. Let $n\in\{3,4\}$ and $k\ge2$. By  Proposition~\ref{lower}, it suffices to  show that $ GR(G_{i_1},  \ldots, G_{i_k}) \le |G_{i_1}|+\sum_{j=2}^k i_j$. \medskip

 By Theorem~\ref{cycles},  Theorem~\ref{path-cycle} and Theorem~\ref{paths}, $ GR(G_{i_1}, G_{i_2})=R(G_{i_1}, G_{i_2})=|G_{i_1}|+i_2$. We may assume that $k\ge3$. 
Let $N: = |G_{i_1}|+\sum_{j=2}^k i_j$. Since $GR_k(P_3)= 3$,  we may assume that $i_1\ge1$ and so  $N\ge2i_1+3\ge5$.  Let  $G$ be a complete graph on $N$ vertices and let $c: E(G)\rightarrow [k]$ be any Gallai  coloring of $ G$ such that all the edges of $G$ are colored by at least three colors under $c$.  We next show that $G$ contains a monochromatic copy of $G_{i_j}$ in color $j$ for some $j\in[k]$. Suppose   $G$ contains no monochromatic copy of $G_{i_j}$ in color $j$ for any $ j\in[k]$ under $c$. Such a Gallai $k$-coloring $c$   is called a  \dfn{bad coloring}.  Among all complete graphs on $ N$ vertices with a   bad coloring,  we choose $G$ with $N$ minimum. \medskip

Consider a Gallai partition of $G$ with parts $A_1,   \ldots, A_{p}$, where  $p\ge2$. 
We may assume that  $|A_1| \ge \cdots \ge |A_p| \ge 1$. 
Let $\mc{R}$ be the reduced graph of $G$ with vertices $a_1,   \ldots, a_p$, where $a_i \in A_i$ for all $i\in[p]$. By Theorem \ref{Gallai}, we may assume that every  edge  of $\mc{R}$
is  colored either  red or  blue.   Since all the edges of $G$ are colored by at least three colors under $c$, we see that $\mc{R}\ne G$ and so $|A_1|\ge2$.   By abusing the notation, we use $i_b$ to denote $i_j$ when the color $j$ is blue. Similarly, we use $i_r$  to denote $i_j$ when the color $j$ is red.
Let
\[
\begin{split}
A_r &:= \{a_j \in \{a_2, \ldots, a_p\} \mid a_ja_1 \text{ is colored red in } \mc{R} \}  \text{ and}\\
A_b &:=  \{a_i \in \{a_2, \ldots, a_p\} \mid a_ia_1 \text{ is colored blue in } \mc{R}\}.
\end{split}
\]
 Let   $R:=\bigcup_{a_j \in A_r} A_j$ and $B:= \bigcup_{a_i \in A_b} A_i$.    
Then   $|A_1|+|R|+|B|=|G|= N$  and $\max\{|B|, |R|\} \ne0$ because $p\ge2$.   Thus $G$ contains a blue $P_3$ between $B$ and $A_1$ or a red $P_3$ between $R$ and $A_1$, and so $\max\{i_b, i_r\} \ge 1$.  We next prove several claims.  \\

\setcounter{counter}{0}

\noindent {\bf Claim\refstepcounter{counter}\label{Observation}  \arabic{counter}.}  Let $r\in[k]$ and let  $s_1, \ldots,s_r$  be nonnegative integers with $s_1+ \cdots+s_r\ge1$.  If $i_{j_1}\geq s_1,  \dots, i_{j_r}\geq s_r$   for colors  $ j_1, j_2, \dots, j_r\in[k]$,   then for any $S\subseteq V(G)$ with $|S|\ge N-(s_1+ \cdots+s_r)$, $G[S]$ must contain a monochromatic copy of $G_{i^*_{j_q}}$ in color $j_q$ for some $j_q\in\{j_1, \ldots,j_r\}$, where   $i^*_{j_q}=i_{j_q}-s_q$.

\pf  Let $i^*_{j_1} :=i_{j_1}-s_1,  \dots,i^*_{j_r} :=i_{j_r}-s_r$,  and $i^*_j :=i_j$ for all $j\in [k]\less\{j_1, \ldots,j_r\}$.  Let   $i^*_\ell:=\max \{i^*_j: {j \in [k]} \} $. Then $i^*_\ell\le i_1$.  Let $N^* :=|G_{i^*_\ell}|+[(\sum_{   j=1}^k i^*_j) -i^*_\ell]$.   Then $N^*\ge 3$ and $ N^*\le N-(s_1+ \cdots+s_r)<N$ because $s_1+\cdots+s_r\ge1$. Since $|S|\ge  N-(s_1+ \cdots+s_r)\ge  N^*$ and  $G[S]$ does not have a monochromatic copy of $G_{i_j}$ in color $j$ for all $j\in [k]\less\{j_1, \ldots,j_r\}$ under $c$,  by   minimality of $N$, $G[S]$ must contain a monochromatic copy of $G_{i^*_{j_q}}$ in color $j_q$ for some $j_q\in\{j_1,\ldots,j_r\}$.\proofsquare

\noindent {\bf Claim\refstepcounter{counter}\label{e:Ap}  \arabic{counter}.}   $|A_1|\le n-1$ and so  $G$ does not contain a  monochromatic copy of  a graph on $|A_1|+1\le n$ vertices in any color $m\in[k]$   that is neither red nor blue.

\pf  Suppose $|A_1| \ge n$. We first claim that $i_b\ge |B|$ and $i_r\ge |R|$. Suppose $i_b\le |B|-1$ or  $i_r\le |R|-1$. Then we obtain  a blue  $G_{i_b}$ using the edges between $B$ and $A_1$  or a red $G_{i_r}$ using  the edges  between $R$ and $A_1$, a contradiction.
Thus $i_b\ge |B|$ and $i_r\ge |R|$,  as claimed.   Let $i_b^*:=i_b-|B|$ and  $i_r^*:=i_r-|R|$. 
Since $|A_1|= N-|B|-|R|$, by Claim~\ref{Observation} applied to $i_b\ge |B|$, $i_r\ge |R|$ and $A_1$, $G[A_1]$ must have a blue $G_{i^*_b}$ or a red   $G_{i^*_r}$, say the latter.  Then  $i_r>i_r^*$. Thus    $|R|>0$ and $G_{i^*_r}$ is a red path on $2i_r^*+3$ vertices. 
 Note that 
 \begin{align*}
|A_1| &=|G_{i_1}| +\sum_{j=2}^k i_j-|B|-|R|\\
&\ge \begin{cases}
|G_{i_r}|+i_b-|B|-|R|   & \ \text{if } \  i_r\ge i_b \\
|G_{i_b}|+i_r-|B|-|R|   & \ \text{if } \  i_r< i_b,  
\end{cases}\\
&\ge
 \begin{cases}
|G_{i_r}|+i^*_b-|R|   & \ \text{if } \  i_r\ge i_b \\
  2i_b+2+i_r-|B|-|R|\ge i^*_b+(2i_r+3)-|R|  & \ \text{if } \  i_r< i_b, 
\end{cases}\\
&\ge  
  |G_{i_r}|-|R|.
\end{align*}
 Then 
\begin{align*}
|A_1|-|G_{i_{r}^*}| & \ge |G_{i_r}|-|G_{i_{r}^*}|-|R|\\
&= 
\begin{cases}
(3+2i_r)-(3+2i_{r}^*)-|R| =|R| & \ \text{if} \  i_r\le n-2 \\
(2+2i_r)-(3+2i_{r}^*)-|R| =|R|-1& \ \text{if} \  i_r= n-1.
\end{cases}
\end{align*}
But then $G[A_1\cup R]$ contains   a red $G_{i_r}$ using the edges of the $G_{i_r^*}$ and the edges between $A_1 \less V(G_{i_r^*})$ and $R$, a contradiction.  This proves that $|A_1|\le n-1$. Next,
let  $m\in[k]$ be  any  color   that is neither red nor blue. Suppose $G$ contains a  monochromatic copy of  a graph, say $J$,  on $|A_1|+1$ vertices in color $m$. Then $V(J)\subseteq A_\ell$ for some $\ell\in[p]$. But then $|A_\ell|\ge|A_1|+1$, contrary to $|A_1|\ge |A_\ell|$. \proofsquare

 For  two disjoint sets $U, W\subseteq V(G)$,  we say $U$  is \dfn{blue-complete} (resp.  \dfn{red-complete})   to $W$    if all the edges between $U$ and $W$    are colored  blue (resp.  red) under $c$.  For convenience, we say $u$  is \dfn{blue-complete} (resp.  \dfn{red-complete})   to $W$  when $U=\{u\}$.\\

\noindent {\bf Claim\refstepcounter{counter}\label{e:R}  \arabic{counter}.}   $\min\{|B|, |R|\} \ge 1$, $p\ge3$  and $B$ is neither red- nor blue-complete to $R$ under $c$.

\pf  Suppose $B=\emptyset$ or $R=\emptyset$. By symmetry, we may assume that $R=\emptyset$.   Then $B\ne \es$ and so $i_b\ge1$. By Claim~\ref{e:Ap}, $|A_1|\le n-1\le 3$ because $n\in\{3,4\}$. Then $|A_1|\le i_b+2$.    If $i_b \le |A_1|-1$, then $i_b\le n-2$ by Claim~\ref{e:Ap}. Thus  $  G_{i_b}$ is a blue path on $2i_b+3$ and so 
 \begin{align*}
 |B|= N-|A_1|\ge |G_{i_b}|-|A_1|=\begin{cases}
i_b+1  & \ \text{if} \ |A_1|= i_b+2\\
i_b+2  & \ \text{if} \ |A_1|= i_b+1.
\end{cases}
\end{align*} 
 But then we obtain a blue   $G_{i_b}$  using the edges between $B$ and $A_1$. Thus    $i_b \ge |A_1|$. Let $i^*_b :=i_b-|A_1|$. By Claim \ref{Observation} applied to $i_b \ge |A_1|$ and $B$, $G[B]$ must have a blue $G_{i^*_b}$. Since  
 \begin{align*}
|B|-|G_{i_{b}^*}| & \ge |G_{i_{b}}|- |G_{i_{b}^*}| -|A_1|= 
\begin{cases}
(3+2i_b)-(3+2i_{b}^*)-|A_1| =|A_1| & \ \text{if} \  i_b\le n-2 \\
(2+2i_b)-(3+2i_{b}^*)-|A_1| =|A_1|-1& \ \text{if} \  i_b= n-1,
\end{cases}
\end{align*}
we see that     $G$ contains  a blue   $G_{i_b}$ using the edges of the $G_{i^*_b}$ and the edges between $B \less V(G_{i^*_b})$ and $A_1$, a contradiction. 
Hence $R\ne \emptyset$ and so $p\ge3$ for any Gallai partition of $G$.   It follows that  $B$ is neither red- nor blue-complete to $R$, otherwise     $\{B\cup A_1, R\}$ or $\{B, R\cup A_1\}$ yields  a  Gallai partition of $G$ with only two parts.
 \proofsquare

 \noindent {\bf Claim\refstepcounter{counter}\label{e:ij}  \arabic{counter}.}   Let  $m\in[k]$ be  a  color   that is neither red nor blue. Then $i_m\le1$. In particular, if $i_m=1$, then $n=4$ and $G$ contains a monochromatic copy of $P_3$  in color $m$ under $c$. 
 
 \pf   By   Claim~\ref{e:Ap},          $G$ contains no   monochromatic copy of $P_n$  in color $m$ under $c$. 
Suppose $i_m\ge1$.    Let $ i^*_m:=i_m-1$. 
 By Claim \ref{Observation} applied to $i_m\ge1$ and  $ V(G) $,  $G$ must have a monochromatic copy of $G_{i^*_m}$  in color $m$ under $c$. Since    $n\in\{3,4\}$ and $G$ contains no   monochromatic copy of $P_n$  in color $m$, we see that $n=4$ and $i^*_m=0$. Thus $i_m=1$ and    $G$ contains  a monochromatic copy of $P_3$  in color $m$ under $c$.\proofsquare

By Claim~\ref{e:R}, $ B\ne \es$ and $R\ne \es$.  Since $|A_1|\ge2$,   we see that $G$ has a blue  $P_3$ using edges between $B$ and $A_1$, and a red $P_3$ using edges between $R$ and $A_1$. Thus $ i_b\ge 1$ and $ i_r\ge1$.   Then $|G_{i_1}|\ge5$ and so   $N= |G_{i_1}|+ \sum_{j=2}^k i_j \ge 6$. By Claim~\ref{e:Ap}, $|A_1|\le n-1$. If   $|B|=|R|=1$,  then $N=|A_1|+|B|+|R|\le n+1 \le 5$, a contradiction. Thus $ |B|\ge2$ or $|R|\ge2$.  
Since  $B$ is neither  red- nor blue-complete to $R$, we see that     $G$ contains either a blue   $P_5$  or a red $P_5$.  Thus $i_1\ge \max\{i_b, i_r\}\ge2\ge n-2$ because $n\in\{3,4\}$.     
By Claim~\ref{e:ij}, we may assume that $\{i_b, i_r\}=\{i_1, i_2\}$. Then 
 \begin{align*}
|G_{i_1}| = 
\begin{cases}
2i_1+2 =1+n+i_1  & \ \text{if} \  i_1=n-1 \\
 2i_1+3 = 1+n+i_1& \ \text{if} \  i_1= n-2.\end{cases}
\end{align*}
  Therefore $N=| G_{i_1}|+\sum_{j=2}^k i_j=1+n+\sum_{j=1}^k i_j\ge 1+n +i_b+i_r$. \\

 \noindent {\bf Claim\refstepcounter{counter}\label{e:R4}  \arabic{counter}.} $ |B|\le n-1$ or  $|R| \le n-1$.  
  
\pf Suppose $ |B|\ge n $ and   $|R| \ge n$.  Let $H=(B,R)$ be the complete bipartite graph obtained from $G[B\cup R]$ by deleting all the edges with both ends in $B$ or both ends in $R$. Then $H$ has no blue  $P_{2n-3}$ with both ends in $B$,  else, we obtain a blue  $C_{2n}$ because $|A_1|\ge2$. Similarly, $H$ has no   red $P_{2n-3}$ with both ends in $R$. For every vertex $v\in B\cup R$,  let $d_b(v): =|\{u: uv \text{ is colored blue in } H\}|$ and $d_r(v):=|\{u: uv \text{ is colored red in } H\}|$.  
Let $x_1, \ldots, x_n \in B$,  $y_1, \ldots, y_n\in R$ and   $a_1, a^*_1\in A_1 $ be all distinct.   We next claim that $d_r(v)\le n-2$ for all $v\in B$. Suppose, say,      $d_r(x_1)\ge n-1$. Then $n=4$ because   $H$ has no red $P_{2n-3}$ with both ends in $R$. 
We may assume that  $x_1$ is red-complete to $\{y_1,  y_2,  y_3\}$. Since $H$ has no red $P_{5}$ with both ends in $R$, we see that for all $i\in\{2,3,4\}$ and every $W\subseteq \{y_1, y_2, y_3\}$ with $|W|=2$, no $x_i$ is red-complete to  $W$. We may further assume that $x_2y_1, x_2y_2, x_3y_1$ are colored blue.  Then $x_4y_2$ must be  colored red, else,  $H$ has a blue $P_5$ with vertices $x_3, y_1, x_2, y_2, x_4$ in order. Thus  $x_4y_1, x_4y_3$ are colored blue.  But then $H$ has a blue $P_5$  with vertices $x_2, y_2, x_3, y_1, x_4$ in order (when $x_3y_2$ is colored blue) or vertices $x_2, y_1, x_3, y_3, x_4$ in order (when $x_3y_3$ is colored blue), a contradiction. Thus $d_r(v)\le n-2$ for all $v\in B$. Similarly, 
$d_b(u)\le n-2$ for all $u\in R$.  
Then $|B||R|=|E(H)|=\sum_{v\in B}d_r(v)+\sum_{u\in R}d_b(u)\le (n-2)|B|+(n-2)|R|$.  Using inequality of arithmetic and geometric means, we obtain  that   $n=4$, $|B|=|R|=4$ and $d_r(v)=d_b(v)=2$ for each $v\in B\cup R$.  Thus the set of  all the blue edges in $H$ induces a $2$-regular spanning subgraph of $H$. Since $H$ has no blue $C_8$, we see that $H$ must contain two vertex-disjoint copies of blue $C_4$. We may assume that $ y_1$ is blue-complete to $\{x_1, x_2\}$ and $ y_2$ is blue-complete to $\{x_3, x_4\}$. But then  $G$ contains   a blue $C_8$ with vertices $a_1, x_1, y_1, x_2, a^*_1, x_3, y_2, x_4$ in order, 
 a contradiction. 
\proofsquare

 \noindent {\bf Claim\refstepcounter{counter}\label{e:A1}  \arabic{counter}.}  $|A_1|=3$ and $n=4$.  

 \pf By   Claim~\ref{e:Ap},    $|A_1|\le n-1\le3$ because $n\in\{3,4\}$. Note that  $|A_1|=3$ only when $n=4$. Suppose $|A_1| = 2$.  By   Claim~\ref{e:Ap},     $G$ has no monochromatic copy of $P_3$   in color $j$ for any $ j\in\{3, \ldots, k\}$ under $c$. By Claim~\ref{e:ij}, $i_3=\cdots=i_k=0$ and so  $N =1+n+\sum_{j=1}^k i_j=1+n +i_b+i_r$. 
We may assume that  $A_1, \ldots, A_t$ are all the parts of order two in the Gallai partition $A_1,  \ldots, A_p$ of $G$, where   $t\in[p]$.  Let $A_i:=\{a_i, b_i\}$ for all $i\in [t]$.  By reordering if necessary,  each of $A_1, \ldots, A_t$   can be chosen as the largest part in the Gallai partition $A_1,    \ldots, A_p$ of $G$. 
For all $i\in [t]$, let
\[
\begin{split}
A^i_b &:=  \{a_j \in V(\mc{R}) \mid a_ja_i \text{ is colored blue in } \mc{R} \} \text{ and }\\
A^i_r &:= \{a_j \in V(\mc{R}) \mid a_ja_i \text{ is colored red in } \mc{R} \}.
\end{split}
\]
 Let $B^i:= \bigcup_{a_j \in A^i_b} A_j$ and $R^i:=\bigcup_{a_j \in A^i_r} A_j$. Then $|B^i|+|R^i|=N-|A_1|=n+i_b+i_r-1\ge n+2$,  because $\max\{i_b, i_r\}\ge2$ and $\min\{i_b, i_r\}\ge1$. Since each of $A_1, \ldots, A_t$   can be chosen as the largest part in the Gallai partition $A_1,   \ldots, A_p$ of $G$,   by Claim \ref{e:R4}, either $|B^i| \le n-1$ or $|R^i| \le n-1$ for all $i\in [t]$. We claim that $|B^i| \ne |R^i|$  for all $i\in [t]$. Suppose $|B^i| = |R^i|$  for some $i\in [t]$.  By Claim \ref{e:R4},  $n+2\le |B^i|+|R^i| \le 2(n-1)\le 6$. It follows that   $|B^i| = |R^i|=3$ and $n=4$. Thus   $G$ has  a blue $P_5$ between $B^i$ and $A_i$  and  a red    $P_5$ between $R^i$ and $A_i$.   It follows that $\min\{i_b, i_r\}\ge2$. But then $|B^i|+|R^i| =n+i_b+i_r-1\ge 7$, a contradiction. This proves that $|B^i| \ne |R^i|$  for all $i\in [t]$. Let
 \[
\begin{split}
E_B :=  \{a_ib_i \mid i\in [t]   \text{ and }  |R^i|<|B^i| \}  \text{ and } E_R :=  \{a_ib_i \mid i\in [t]   \text{ and }  |R^i|>|B^i| \}.  
\end{split}
\]
We next apply the recoloring method. Let $c^*$ be  an edge-coloring of $G$  obtained from $c$ by recoloring all the edges in $E_B$   blue and all the edges in $E_R$   red. Then every edge of  $G$ is  colored either red or blue under $c^*$. Since $|G|=1+n +i_b+i_r\ge R(G_{i_b}, G_{i_r})$ by Theorem~\ref{cycles}, Theorem~\ref{path-cycle} and Theorem~\ref{paths}, we see that $G$ must contain a blue  $G_{i_b}$ or a red $G_{i_r}$ under $c^*$.  By symmetry, we may   assume that $G$ has a blue $H:=G_{i_b}$ under $c^*$.  Then $H$ contains no edges   of $E_R$ but must contain  at least one edge of $E_B$, else,  we obtain a blue $G_{i_b}$ in $G$ under $c$.  We choose $H$ so that $|E(H) \cap E_B|$ is minimal. We may further assume that $a_1b_1 \in E(H)$. By the choice of $c^*$, $|R^1|\le n-1$ and $|R^1|<|B^1|$. Then $|B^1| \ge 2$ and so $G$ has a blue $P_5$ under $c$  because $B^1$ is not red-complete to $R^1$. Thus $i_b \ge 2$.   Let  $W:=V(G) \less V(H) $.  \medskip

We next claim that  $i_b=n-1$. Suppose  $2 \le i_b\le n-2$. Then $n=4$, $i_b=2$, $H=P_7$ and $|G| = 1+n+i_b +i_r = 7+i_r $. Thus  $|W|= i_r$. Let $x_1,  \ldots, x_7$ be the vertices of $H$ in order.  By symmetry, we may assume that $x_\ell x_{\ell+1}=a_1b_1$ for some $\ell\in [3]$. Then $W\cup\{x_7\}$ must be  red-complete to $\{a_1, b_1\}$ under $c$, else,  say a vertex $u\in W\cup\{x_7\}$, is blue-complete to $\{a_1, b_1\}$ under $c$, then we obtain a blue $H':=P_7$ under $c^*$ with  vertices $  x_1, \ldots, x_{\ell}, u,  x_{\ell+1}, \ldots, x_6$  in order    such that $|E(H') \cap E_B| < |E(H) \cap E_B|$,  contrary to the choice of $H$. 
Thus $W\cup\{x_7\}\subseteq R^1$ and so  $|R^1| \ge |W\cup\{x_7\}|=i_r+1 \ge 2$. Note that   $G$ contains a red $P_5$ under $c$  because $|R^1|\ge2$ and $R^1$ is not blue-complete to $B^1$. Thus $i_r \ge 2$. 
Then $3\le i_r+1 \le |R^1| \le 3$, which implies that  $i_r=2$ and $R^1=W\cup\{x_7\}$. Thus  $\{a_1, b_1\}$ is blue-complete to $V(H) \less \{x_\ell, x_{\ell+1}, x_7\}$. But then we
obtain a blue $H':=P_7$ under $c^*$  with  vertices $ x_1, \ldots, x_{\ell},  x_{\ell+2}, x_{\ell+1}, x_{\ell+3}, \ldots,  x_7$  in order   such that $|E(H') \cap E_B| < |E(H) \cap E_B|$, a contradiction.  
This proves that $i_b=n-1$.  \medskip

Since $i_b=n-1$, we see that $H=C_{2n}$.  Then $|G|=1+n +i_b+i_r=2n+i_r$ and so $|W|=i_r$.
   Let $a_1, x_1, \ldots,  x_{2n-2}, b_1$ be the vertices of $H$ in order and let $W=V(G) \less V(H):=\{w_1, \ldots,  w_{i_r}\}$. Then  $x_1b_1$ and $a_1x_{2n-2}$ are colored blue under $c$ because $\{a_1, b_1\}=A_1$.  Suppose  $\{x_j, x_{j+1}\}$ is   blue-complete to $\{a_1, b_1\}$ under $c$ for some $j\in[2n-3]$. Then $G$ has   a blue $H':=C_{2n}$ under $c^*$ with  vertices $a_1, x_1, \ldots, x_j, b_1,  x_{2n-2}, \ldots, x_{j+1}$ in order   such that $|E(H') \cap E_B| < |E(H) \cap E_B|$, contrary to the choice of $H$. Thus, for all $j\in[2n-3]$, $\{x_j, x_{j+1}\}$ is not blue-complete to $\{a_1, b_1\}$. Since $\{x_1, x_{2n-2}\}$ is blue-complete to $\{a_1, b_1\}$ under $c$, we see    that $x_2, x_{2n-3}\in R^1$ and then $|R^1\cap \{x_2, \ldots, x_{2n-3}\}|= |R^1|=n-1$. Thus $R^1=\{x_2, x_3\}$ when $n=3$. By symmetry, we may assume that  $R^1=\{x_2, x_3, x_{5}\}$  when $n=4$.  Then  $W \subseteq B^1$. Thus  $R^1$ is red-complete to $\{a_1, b_1\}$   and $W$ is blue-complete to $\{a_1, b_1\}$ under $c$.      It follows   that     for any $w_j\in W$ and $x_m\in R^1$, $\{x_m, w_j\}\ne A_i$ for all  $i\in [t]$.     
Then $x_2$  must be   red-complete to $W$ under $c$, else, say $x_2w_1$ is colored blue  under $c$, then we obtain a blue $H':=C_{2n}$  under $c^*$ with vertices $a_1, x_1, x_2, w_1,  b_1,  x_4$ (when $n=3$) and vertices $a_1, x_1, x_2, w_1,  b_1,  x_4, x_5, x_6$ (when $n=4$)  in order  such that $|E(H') \cap E_B| < |E(H) \cap E_B|$,  a contradiction.    Similarly,  $x_{3}$ is red-complete to $W$ under $c$, else, say $x_3w_1$ is colored blue  under $c$, then we obtain a blue $H':=C_{2n}$ under $c^*$ with vertices $b_1, x_4, x_3, w_1,  a_1,  x_1$  (when $n=3$) and vertices $b_1, x_6, x_5, x_4, x_3, w_1,  a_1,  x_1$ (when $n=4$)  in order  such that $|E(H') \cap E_B| < |E(H) \cap E_B|$,  a contradiction. Thus $\{x_2, x_{3} \}$ is red-complete to $W$ under $c$. 
Then  for any $w_j\in W$, $\{x_1, w_j\}\ne A_i$ for all  $i\in [t]$ since   $x_2x_1$ is colored blue and $x_2$ is red-complete to $W$ under $c$.  If $x_1w_j$ is colored blue under $c$ for some $w_j\in W$, then we obtain a blue $H':=C_{2n}$ under $c^*$ with vertices $a_1, w_j,  x_1,  \ldots, x_{2n-2}$ in order  such that $|E(H') \cap E_B| < |E(H) \cap E_B|$,  a contradiction. Thus   $\{x_1, x_2, x_{3} \}$ is red-complete to $W$ under $c$. Then $|W|= i_r\ge2$  because $G$ contains a red $P_5$ under $c$ with vertices $x_1, w_1, x_2, a_1, x_{3} $ in order. 
But then we obtain a red $ C_{2n}$ under $c$ with vertices  $a_1, x_2, w_1, x_1, w_2, x_3$ in order (when $n=3$) and $a_1, x_2, w_1, x_1, w_2, x_3, b_1, x_5$ in order (when $n=4$), a contradiction. \proofsquare

By Claim~\ref{e:A1}, $|A_1| =3$ and $n=4$. Then $|B\cup R|=N-|A_1| \ge 2+i_b+i_r \ge 5$ because $\max\{i_b, i_r\}\ge2$ and $\min\{i_b, i_r\}\ge1$. By symmetry, we may assume that $|B| \ge |R|$. Then $|B| \ge 3$ and so  $G$ has a blue $P_7$ because $|A_1| =3$ and $B$ is not red-complete to $R$. Thus  $i_b=3$. By Claim~\ref{e:R4}, $|R| \le 3$. Then $i_r \ge |R|$, else,  we obtain a red $G_{i_r}$ because $|A_1| =3$ and $R$ is not blue-complete to $B$. 
Then $|B|\ge 2+i_b+i_r  -|R| \ge 5$.   Thus  $G[B \cup R]$ has no blue $P_3$ with both ends in $B$, else, we obtain a blue $C_8$ because $|A_1|=3$ and $|B|\ge5$. %
Let $i^*_b:=0$ and $i^*_r:=i_r-|R|\le2$. By Claim~\ref{Observation} applied to $i_b = |A_1| $, $i_r \ge |R|$ and $B$, $G[B]$ must contain a red $ P_{2i^*_r+3}$ with vertices, say $x_1,  \ldots, x_{2i^*_r+3}$, in order.  Let $  R:=\{y_1, \ldots, y_{|R|}\}$. Then no $y_j\in R$ is   blue-complete to any $W\subseteq B$ with $|W|=2$, in particular, when $W=\{x_1, x_{2i^*_r+3}\}$, because $G[B \cup R]$ has no blue $P_3$ with both ends in $B$. We may assume that $x_1y_1$ is colored red. Note that $G[R\cup A_1]$ has a red $P_{2|R|}$ with $y_1$ as an end.  Then $G[\{x_1,  \ldots, x_{2i^*_r+3}\}\cup R\cup A_1]$ has a red $P_{2i_r+3}$.  It follows that $i_r=3$.  Let $ a^*_1\in A_1 \less\{a_1\}$.  \medskip

Suppose first that  $x_{2i^*_r+3}$ is blue-complete to $R=\{y_1, \ldots, y_{|R|}\}$. Since $G[B \cup R]$ has no blue $P_3$ with both ends in $B$, we see that   $\{x_{2i^*_r+3}\}=A_\ell$ for some $\ell\in[p]$,      $B \less \{x_{2i^*_r+3}\}$ is red-complete to $ \{y_1, \ldots, y_{|R|}\}$, and $x_{2i^*_r+3}$ is adjacent to at most one vertex, say $w\in B$, such that $wx_{2i^*_r+3}$ is colored blue. 
Thus   $x_{2i^*_r+3}$  is red-complete to $B\less\{w, x_{2i^*_r+3}\}$. Let    $w^*\in B\less\{x_1, x_2, x_3, w\}$. Since $B \less \{x_{2i^*_r+3}\}$ is red-complete to $ \{y_1, \ldots, y_{|R|}\}$, we see that $\{x_1, \ldots, x_{2i^*_r+2}\}$ is red-complete to $ \{y_1, \ldots, y_{|R|}\}$.  If $w\notin \{x_2, \ldots, x_{2i^*_r+1}\}$, then we obtain a red $C_8$ with vertices $y_1, x_1, x_2, x_7, x_3,   \ldots, x_{6}$ (when $i^*_r=2$),     vertices $a_1, y_1, x_1, x_2, x_5, x_3, x_4,  y_2$ (when $i^*_r=1$), and vertices $a_1, y_1, x_2, x_3, w^*,  y_2, a^*_1, y_3$ (when $i^*_r=0$) in order, a contradiction. Thus $w\in \{x_2, \ldots, x_{2i^*_r+1}\}$. Then  $i^*_r\ge1$ and $x_1x_{2i^*_r+1}$ is colored red.  But then we obtain a red $C_8$ with vertices $  y_1, x_2, x_3, x_4, x_5, x_6,x_7, x_1 $ (when $i^*_r=2$) and  vertices $a_1, y_1, x_2, x_3, x_4, x_5, x_{1},  y_2$ (when $i^*_r=1$) in order, a contradiction. 
This proves that  $x_{2i^*_r+3}$ is not blue-complete to $R$. Then $|R|\ge2$, else, $|R|=1$, $i^*_r=2$ and $x_7y_1$ is colored red, which yields a red $C_8$ with vertices $y_1, x_1, \ldots, x_7$ in order, a contradiction. 
Thus $i^*_r\le1$. Next, suppose $x_{2i^*_r+3}$ is not blue-complete to $\{y_2, \ldots, y_{|R|}\}$, say $x_{2i^*_r+3}y_2$ is colored red. By assumption, $x_1y_1$ is red. We then obtain a red $C_8$ with vertices $a_1, y_1, x_1,   \ldots, x_{5}, y_2$ (when $i^*_r=1$) and    vertices $a_1, y_1, x_1,  x_2, x_3,   y_2, a^*_1, y_3$ (when $i^*_r=0$) in order, a contradiction. Thus  $ x_{2i^*_r+3}$ is   blue-complete to $\{y_2, \ldots, y_{|R|}\}$ and so $ x_{2i^*_r+3}y_1$ is colored red. By symmetry of $x_1$ and $ x_{2i^*_r+3}$, $x_{1}$ must be   blue-complete to  $\{y_2, \ldots, y_{|R|}\}$. But then $G[B\cup R]$ has  a blue $P_3$ with vertices $x_1, y_2, x_{2i^*_r+3}$ in order, a contradiction.\medskip

 This completes the proof of Theorem~\ref{main}.\proofsquare

\section*{Acknowledgements}
The authors would like to thank Christian Bosse for  many helpful comments and discussion. We
also thank referees for their careful reading and many helpful comments. In particular,
we are indebted to one referee for an improved version of the formula given in Proposition 1.9, which greatly improved the proof of Theorem \ref{main}.

\end{document}